\documentclass{elsart}
\usepackage{amsfonts}
\usepackage{latexsym}
\usepackage{amsmath}
\usepackage{amssymb}

\newcommand{\C}{\mathbb{C}}

\newcommand{\Z}{\mathbb{Z}}

\newtheorem{defin}{Definition}[section]

\newtheorem{theorem}[defin]{Theorem}
\newtheorem{exa}[defin]{Example}

\newtheorem{lemma}[defin]{Lemma}

\newtheorem{corollary}[defin]{Corollary}
\newenvironment{proof}
{\noindent{\it Proof.}}{\hfill $\Box$\par\vspace{2.5mm}}

\begin{document}

\begin{frontmatter}

\title{Difference analogue of the Lemma on the Logarithmic Derivative with
applications to difference equations\thanksref{label1}}
\thanks[label1]{The research reported in this paper was supported in
part by EPSRC grant number GR/R92141/01 and by the Finnish Academy
grant number 204819.}

\author{R.G. Halburd},
\address{
Department of Mathematical Sciences,
Loughborough University\\
Loughborough, Leicestershire, LE11 3TU, UK
\\r.g.halburd@lboro.ac.uk
}

\author{R.J. Korhonen\corauthref{cor}}
\corauth[cor]{Corresponding author.}
\address{
Department Mathematics,
University of Joensuu\\
P.O. Box 111, FI-80101 Joensuu, Finland\\risto.korhonen@joensuu.fi
}

\begin{abstract}
The Lemma on the Logarithmic Derivative of a meromorphic function
has many applications in the study of meromorphic functions and
ordinary differential equations. In this paper, a difference
analogue of the Logarithmic Derivative Lemma is presented, and
then applied to prove a number of results on meromorphic solutions
of complex difference equations. These results include a
difference analogue of the Clunie Lemma, as well as other results
on the value distribution of solutions.
\end{abstract}


\end{frontmatter}

\section{Introduction}

The Lemma on the Logarithmic Derivative states that outside of a
possible small exceptional set
    \begin{equation}\label{logder}
    m\left(r,\frac{f'}{f}\right)=O(\log T(r,f)+\log r),
    \end{equation}
where $m(r,f)$ denotes the Nevanlinna proximity function and
$T(r,f)$ is the characteristic of a meromorphic function $f$
\cite{hayman:64}. This is undoubtedly one of the most useful
results of Nevanlinna theory, having a vast number of applications
in the theory of meromorphic functions and in the theory of
ordinary differential equations. For instance, the proofs of the
Second Main Theorem of Nevanlinna theory \cite{nevanlinna:25} and
Yosida's generalization \cite{yosida:33} of the Malmquist theorem
\cite{malmquist:13} both rely heavily on the Lemma on the
Logarithmic Derivative. One major problem in the study of complex
difference equations has so far been the lack of efficient tools,
which can play roles similar to that played by relation
\eqref{logder} for differential equations. This has meant that
most results have had to be proved separately for each difference
equation. This slows down the efforts to construct a coherent
theory, and it may be one of the reasons why the theory of
meromorphic solutions of complex difference equations is not as
developed as the theory of differential equations.

The foundations of the theory of complex difference equations were
laid by N\"orlund, Julia, Birkhoff, Batchelder and others in the
early part of the twentieth century. Later on,
Shimomura~\cite{shimomura:81} and
Yanagihara~\cite{yanagihara:80,yanagihara:85} studied non-linear
complex difference equations from the viewpoint of Nevanlinna
theory. Recently, there has been a renewed interested in the
complex analytic properties of solutions of difference equations.
In differential equations, Painlev\'e and his colleagues
identified all equations, out of a large class of second order
ordinary differential equations, that possess the Painlev\'e
property \cite{fuchs:05,gambier:10,painleve:02}. Those equations
which could not be integrated in terms of known functions or
through solutions of linear equations are now known as the
Painlev\'e differential equations. Similarly, Ablowitz, Halburd
and Herbst \cite{ablowitzhh:00} suggested that the growth of
meromorphic solutions of difference equations could be used to
identify those equations which are of ``Painlev\'e type''. In
\cite{halburdk:04} the existence of one finite order non-rational
meromorphic solution was shown to be sufficient to reduce a
general class of second-order difference equations to one of
difference Painlev\'e equations or to a linear difference
equation, provided that the solution does not satisfy a certain
first-order difference Riccati equation. The proof of this fact
relies on a difference analogue of the Lemma on the Logarithmic
Derivative, Theorem~\ref{logdiff} below, as well as on its
consequences, Theorems~\ref{clunieanalogue} and
\ref{mohonkoanalogue}, which were used in \cite{halburdk:04}
without proving them. Findings in \cite{halburdk:04} suggest that
finite order meromorphic solutions of difference equations have a
similar role as meromorphic solutions of differential equations.

The purpose of this paper is to prove a difference analogue of the
Lemma on the Logarithmic Derivative, and to apply it to study
meromorphic solutions of large classes of difference equations.
The difference analogue appears to be in its most useful form when
applied to study finite order meromorphic solutions of difference
equations, which is in agreement with the findings in
\cite{halburdk:04}. Applications include, for instance, a
difference analogue of the Clunie Lemma \cite{clunie:62}. The
original lemma has proved to be an invaluable tool in the study of
non-linear differential equations. The difference analogue gives
similar information about the finite order meromorphic solutions
of non-linear difference equations.

\section{Difference analogue of the Lemma on the Logarithmic Derivative}

\begin{theorem}\label{logdiff}
Let $f$ be a non-constant meromorphic function, $c\in\C$, $\delta<1$ and $\varepsilon>0$. Then
    \begin{equation}\label{logdiffeq}
    m\left(r,\frac{f(z+c)}{f(z)}\right)= o\left(\frac{T(r+|c|,f)^{1+\varepsilon}}{r^\delta}\right)
    \end{equation}
for all $r$ outside of a possible exceptional set $E$ with finite logarithmic measure $\int_E\frac{dr}{r}<\infty$.
\end{theorem}

\begin{proof} Let $\xi(x)$ and $\phi(r)$ be positive, nondecreasing,
continuous functions defined for $e\leq x<\infty$ and $r_0\leq
r<\infty$, respectively, where $r_0$ is such that $T(r+|c|,f)\geq
e$ for all $r\geq r_0$.  Then by Borel's Lemma \cite[Lemma
3.3.1]{cherryy:01}
    \begin{equation*}
    T\left(r+|c|+\frac{\phi(r)}{\xi(T(r+|c|,f))},f\right) \leq 2 T(r+|c|,f)
    \end{equation*}
for all $r$ outside of a set $E$ satisfying
    \begin{equation*}
    \int_{E\cap [r_0,R]} \frac{dr}{\phi(r)} \leq \frac{1}{\xi(e)}+ \frac{1}{\log
     2}\int_e^{T(R+|c|,f)}\frac{dx}{x\xi(x)}
    \end{equation*}
where $R<\infty$. Therefore, by choosing $\phi(r)=r$ and $\xi(x)=x^{\varepsilon/2}$ with $\varepsilon>0$, and defining
    \begin{equation}\label{alpha}
    \alpha = 1+\frac{r}{(r+|c|)T(r+|c|,f)^\frac{\varepsilon}{2}},
    \end{equation}
 we have
    \begin{equation}\label{estimate}
    T(\alpha(r+|c|),f)=T\left(r+|c|+\frac{\phi(r)}{\xi(T(r+|c|,f))},f\right) \leq  2T(r+|c|,f)
    \end{equation}
for all $r$ outside of a set $E$ with finite logarithmic measure.
Hence, if $f(0)\not=0,\infty$, the assertion follows by combining
\eqref{alpha} and \eqref{estimate} with Lemma~\ref{details} below.
Otherwise we apply Lemma~\ref{details} with the function
$g(z)=z^pf (z)$, where $p\in\Z$ is chosen such that
$g(0)\not=0,\infty$.
\end{proof}

When $f$ is of finite order, the right side of
\eqref{logdiffeq} is small compared to $T(r,f)$, and therefore
relation~\eqref{logdiffeq} is a natural analogue of the Lemma on the
Logarithmic Derivative \eqref{logder}. Concerning the sharpness of
Theorem~\ref{logdiff}, the finite order functions $\Gamma(z)$,
$\exp(z^n)$ and $\tan(z)$ show that $\delta$ in \eqref{logdiffeq}
cannot be replaced by any number strictly greater than one.

If $f$ is of infinite order, the quantity $T(r+|c|,f)
r^{-\delta}$ may be comparable to $T(r,f)$. For instance, by
choosing $f(z)=\exp(\exp(z))$, we have
    \begin{equation*}
    m\left(r,\frac{f(z+1)}{f(z)}\right) = (e-1)T(r,f).
    \end{equation*}
Therefore Theorem \ref{logdiff} is mostly useful when applied to
functions with finite order, although the assertion remains valid
for all meromorphic functions. In the finite-order case we can also remove the $\varepsilon$ in Theorem~\ref{logdiff}.

\begin{corollary}\label{logdiffcor}
Let $f$ be a non-constant meromorphic function of finite order, $c\in\C$ and $\delta<1$. Then
    \begin{equation}\label{logdiffeqcor}
    m\left(r,\frac{f(z+c)}{f(z)}\right)= o\left(\frac{T(r+|c|,f)}{r^\delta}\right)
    \end{equation}
for all $r$ outside of a possible exceptional set with finite logarithmic measure.
\end{corollary}

\begin{proof} Choose any $\delta <1$ and denote
$\delta'=(1+\delta)/2$. Since $f$ is of finite order, we have
$T(r+|c|,f)\leq r^{\rho}$ for some $\rho>0$ and for all $r$
sufficiently large. Therefore, by Theorem~\ref{logdiff}
    \begin{equation*}
    m\left(r,\frac{f(z+c)}{f(z)}\right)= o\left(\frac{T(r+|c|,f)}{r^{\delta'-\varepsilon\rho}}\right),
    \end{equation*}
where $\varepsilon>0$. The assertion follows by choosing $\varepsilon=(\delta'-\delta)/\rho$.
\end{proof}

Note that by replacing $z$ by $z+h$, where $h\in\C$, and $c$ by
$c-h$ in \eqref{logdiffeqcor}, and using the inequality
    \begin{equation*}
    T(r,f(z+h)) \leq (1+\varepsilon)T(r+|h|,f(z)), \qquad \varepsilon>0, \qquad r>r_0,
    \end{equation*}
see \cite{yanagihara:80} or \cite{ablowitzhh:00}, we immediately have
    \begin{equation}\label{logdiffeq2}
    m\left(r,\frac{f(z+c)}{f(z+h)}\right)= o\left(\frac{T(r+|c-h|+|h|,f)}{r^\delta}\right)
    \end{equation}
for all $\delta<1$ outside of a possible exceptional set with finite
logarithmic measure.

\begin{lemma}\label{details}
Let $f$ be a meromorphic function such that $f(0)\not=0,\infty$
and let $c\in\C$. Then for all $\alpha >1$, $\delta<1$ and $r\geq1$,
    \begin{equation*}
    m\left(r,\frac{f(z+c)}{f(z)}\right) \leq
    \frac{K(\alpha,\delta,c)}{r^\delta}\left(T\big(\alpha(r+|c|),f\big)+\log^{+}\frac{1}{|f(0)|}\right),
    \end{equation*}
where
    \begin{equation*}
    K(\alpha,\delta,c)= \frac{8|c|(3\alpha+1)+8\alpha(\alpha-1)|c|^\delta}{\delta(1-\delta)(\alpha-1)^2 r^\delta}.
    \end{equation*}
\end{lemma}

\begin{proof} Let $\{a_n\}$ denote the sequence of all zeros of $f$,
and similarly let $\{b_m\}$ be the pole sequence of $f$, where
$\{a_n\}$ and $\{b_m\}$ are listed according to their
multiplicities and ordered by increasing modulus.  By applying
Poisson-Jensen formula with $s=\frac{\alpha+1}{2}(r+|c|)$, see,
for instance, \cite[Theorem 1.1]{hayman:64}, we obtain
    \begin{equation}\label{meq}
    \begin{split}
    \log \left|\frac{f(z+c)}{f(z)}\right| &= \int_0^{2\pi}
    \log|f(se^{i\theta})|\textrm{Re}\left(\frac{se^{i\theta}+z+c}{se^{i\theta}-z-c}-
    \frac{se^{i\theta}+z}{se^{i\theta}-z}\right)\,\frac{d\theta}{2\pi}\\
    &\quad  + \sum_{|a_n|<s} \log \left|\frac{s(z+c-a_n)}{s^2-\bar a_n(z+c)}\frac{s^2-
    \bar    a_n z)}{s(z-a_n)}\right| \\
    &\quad  - \sum_{|b_m|<s} \log \left|\frac{s(z+c-b_m)}{s^2-\bar b_m(z+c)}\frac{s^2-
    \bar    b_m z)}{s(z-b_m)}\right|\\
    & =: S_1(z) + S_2(z) - S_3(z).
    \end{split}
    \end{equation}
Therefore, by denoting
$E:=\{\varphi\in[0,2\pi):\left|\frac{f(re^{i\varphi}+c)}{f(re^{i\varphi})}\right|\geq
1\}$, we have
    \begin{equation*}
    \begin{split}
    m\left(r,\frac{f(z+c)}{f(z)}\right) &= \int_E \log
    \left|\frac{f(re^{i\varphi}+c)}{f(re^{i\varphi})}\right|\, \frac{d\varphi}{2\pi} \\
    & \leq
     \int_0^{2\pi}|S_1(re^{i\varphi})|+|S_2(re^{i\varphi})|+|S_3(re^{i\varphi})|\,\frac{d\varphi}{2\pi}.
     \end{split}
     \end{equation*}
We will now proceed to estimate each $\int_0^{2\pi}
|S_j(re^{i\varphi})|\,\frac{d\varphi}{2\pi}$ separately. Since
    \begin{equation*}
    \begin{split}
    |S_1| &= \left| \int_0^{2\pi}
    \log|f(se^{i\theta})|\textrm{Re}\left(\frac{2c se^{i\theta}}{(se^{i\theta}-z-c)(se^{i\theta}-
    z)}\right)\,\frac{d\theta}{2\pi}\right|\\
    &\leq   \frac{2|c|s}{(s-r-|c|)^2} \int_0^{2\pi}\left|
    \log|f(se^{i\theta})|\right|\,\frac{d\theta}{2\pi}\\
    &= \frac{2|c|s}{(s-r-|c|)^2} \left(m(s,f)+m\left(s,\frac{1}{f}\right)\right),
    \end{split}
    \end{equation*}
we have
     \begin{equation}\label{S1}
    \int_0^{2\pi} |S_1(re^{i\varphi})|\,\frac{d\varphi}{2\pi}\leq \frac{4|c|s}{(s-r-
    |c|)^2}\left(T(s,f)+\log^{+}\frac{1}{|f(0)|}\right).
     \end{equation}

Next we consider the cases $j=2,3$ combined together. First, by
denoting $\{q_k\}:=\{a_n\}\cup\{b_m\}$ and using the fact that
$|\log x| = \log^{+} x + \log^{+}(1/x)$ for all $x>0$, we have
    \begin{equation}\label{S23}
    \begin{split}
    \int_0^{2\pi} |S_2(re^{i\varphi})| & +|S_3(re^{i\varphi})|\,\frac{d\varphi}{2\pi} \leq \sum_{|q_k|<s}
    \int_0^{2\pi}\log^{+}\left|1+\frac{c}{re^{i\theta}-q_k}\right| \,\frac{d\theta}{2\pi}
     \\
    &  +  \sum_{|q_k|<s}\int_0^{2\pi}\log^{+}\left|1-\frac{c}{re^{i\theta}+c-q_k}\right|
    \,\frac{d\theta}{2\pi}
    \\
    & +  \sum_{|q_k|<s} \int_0^{2\pi}\log^{+}\left|1+\frac{\bar q_k c}{s^2-\bar
    q_k(z+c)}\right| \,\frac{d\theta}{2\pi} \\
    &+ \sum_{|q_k|<s}  \int_0^{2\pi}\log^{+}\left|1-\frac{\bar q_k c}{s^2-\bar
    q_k z}\right| \,\frac{d\theta}{2\pi}.
    \end{split}
    \end{equation}
Second, for any $a\in\C$,  and for all $\delta<1$,
    \begin{equation*}
    \int_0^{2\pi}\frac{d\theta}{|re^{i\theta}-a|^\delta} \leq
    4\int_0^{\frac{\pi}{2}}\frac{d\theta}{|re^{i\theta}-|a||^\delta} \leq \frac{2\pi}{1-
    \delta}\frac{1}{r^\delta}
    \end{equation*}
since $|re^{i\theta}-|a|| \geq r\theta\frac{2}{\pi}$ for all
$0\leq\theta\leq\frac{\pi}{2}$. Therefore
    \begin{equation}\label{eq1}
    \begin{split}
    \int_0^{2\pi}\log^{+}\left|1+\frac{c}{re^{i\theta}-a}\right|\,\frac{d\theta}{2\pi} &\leq
    \frac{1}{\delta}\int_0^{2\pi} \log^{+}\left(1+\left|\frac{c}{re^{i\theta}-
    a}\right|^\delta\right) \,\frac{d\theta}{2\pi} \\
    & \leq \frac{1}{\delta}\int_0^{2\pi} \left|\frac{c}{re^{i\theta}-
    a}\right|^\delta \,\frac{d\theta}{2\pi} \\
    &\leq \frac{|c|^\delta}{\delta(1-\delta)}\frac{1}{r^\delta},
    \end{split}
    \end{equation}
and similarly
    \begin{equation}\label{eq2}
    \int_0^{2\pi}\log^{+}\left|1-\frac{c}{re^{i\theta}+c-a}\right|\,\frac{d\theta}{2\pi}\leq
   \frac{|c|^\delta}{\delta(1-\delta)}\frac{1}{r^\delta}.
    \end{equation}
Third, since for all $a$ such that $|a|<s$,
    \begin{equation*}
    \left|\frac{a}{s^2-\bar a z}\right| \leq \frac{1}{s-r},
    \end{equation*}
we have
    \begin{equation}\label{eq3}
    \int_0^{2\pi}\log^{+}\left|1+\frac{\bar a c}{s^2-\bar
    a(z+c)}\right| \,\frac{d\theta}{2\pi} \leq \frac{|c|}{s-r-|c|}
    \end{equation}
and
    \begin{equation}\label{eq4}
    \int_0^{2\pi}\log^{+}\left|1-\frac{\bar a c}{s^2-\bar
    a z}\right| \,\frac{d\theta}{2\pi} \leq \frac{|c|}{s-r}.
    \end{equation}
Finally, by combining inequalities \eqref{S1} -- \eqref{eq4}, we
obtain
    \begin{equation*}
    \begin{split}
     m\left(r,\frac{f(z+c)}{f(z)}\right) & \leq
      \left( \frac{2|c|}{s-
    r-|c|}+\frac{2|c|^\delta}{\delta(1-\delta)}\frac{1}{r^\delta}\right)
    \left(n(s,f)+n\left(s,\frac{1}{f}\right)\right)\\ & \quad +
     \frac{4|c|s}{(s-r-
    |c|)^2}\left(T(s,f)+\log^{+}\frac{1}{|f(0)|}\right).
    \end{split}
    \end{equation*}
Therefore, using the fact that
    \begin{equation*}
    n(s,f)+n\left(s,\frac{1}{f}\right) \leq \frac{4\alpha}{\alpha-
    1}\left(T(\alpha(r+|c|),f)+\log^{+}\frac{1}{|f(0)|}\right),
    \end{equation*}
see \cite[p. 37]{hayman:64}, and $s=\frac{\alpha+1}{2}(r+|c|)$, we conclude
    \begin{equation*}
    \begin{split}
     &m\left(r,\frac{f(z+c)}{f(z)}\right) \\ &\leq \left(\frac{8|c|(3\alpha+1)}{(\alpha-
     1)^2(r+|c|)}+\frac{8\alpha|c|^\delta}{\delta(1-\delta)(\alpha-1)
     r^\delta}\right)\left(T(\alpha(r+|c|),f)+\log^{+}\frac{1}{|f(0)|}\right)\\
     &\leq \frac{8|c|(3\alpha+1)+8\alpha(\alpha-1)|c|^\delta}{\delta(1-\delta)(\alpha-1)^2 r^\delta}
     \left(T(\alpha(r+|c|),f)+\log^{+}\frac{1}{|f(0)|}\right).
     \end{split}
    \end{equation*}
\end{proof}

\section{Difference analogues of the Clunie and Mohon'ko lemmas}\label{diffsec}

The Lemma on the Logarithmic Derivative is an integral part of the
proof of the Second Main Theorem, one of the deepest results of
Nevanlinna theory. In addition,  logarithmic derivative estimates
are crucial for applications to complex differential equations.
Similarly, Theorem~\ref{logdiff} enables an efficient study of
complex analytic properties of finite order meromorphic solutions of
difference equations. We are concerned with functions which are
polynomials in $f(z+c_j)$, where $c_j\in\C$, with coefficients
$a_\lambda(z)$ such that
    \begin{equation*}
    T(r,a_\lambda)=o(T(r,f))
    \end{equation*}
except possibly for a set of $r$ having finite logarithmic measure.
Such functions will be called \textit{difference polynomials in}
$f(z)$. We also denote
    \begin{equation*}
    |c|:=\max\{|c_j|\}.
    \end{equation*}
The following theorem is analogous to the Clunie Lemma
\cite{clunie:62}, which has numerous applications in the study of
complex differential equations, and beyond.

\begin{theorem}\label{clunieanalogue}
Let $f(z)$ be a non-constant meromorphic solution of
    \begin{equation*}
    f(z)^n P(z,f)=Q(z,f),
    \end{equation*}
where $P(z,f)$ and $Q(z,f)$ are difference polynomials in $f(z)$, and let $\delta<1$ and $\varepsilon>0$.
If the degree of $Q(z,f)$ as a polynomial in $f(z)$ and its shifts
is at most $n$, then
    \begin{equation*}
    m\big(r,P(z,f)\big) = o\left(\frac{T(r+|c|,f)^{1+\varepsilon}}{r^\delta}\right) +o(T(r,f))
    \end{equation*}
for all $r$ outside of a possible exceptional set with finite logarithmic measure.
\end{theorem}

\begin{proof} We follow the reasoning behind the original Clunie
Lemma, see, for instance, \cite{hayman:64,laine:93}, just
replacing the Lemma on the Logarithmic Derivative with
Theorem~\ref{logdiff}. First of all,
    \begin{equation}\label{mP}
    m(r,P) = \int_{E_1} \log^{+}|P(re^{i\theta},f)|\,\frac{d\theta}{2\pi} +      \int_{E_2} \log^{+}|P(re^{i\theta},f)|\,\frac{d\theta}{2\pi},
    \end{equation}
where $E_1=\{\theta\in[0,2\pi]:|f(re^{i\theta})|<1\}$, and $E_2$ is the complement of $E_1$. Now, by denoting $P(z,f)=\sum_\lambda a_\lambda(z) F_\lambda(z,f)$, we have
    \begin{equation*}
    |a_\lambda(re^{i\theta})F_\lambda(re^{i\theta},f)| \leq
    |a_\lambda(re^{i\theta})|\left|\frac{f(re^{i\theta}+c_1)}{f(re^{i\theta})}\right|^{l_1}\cdots\left|\frac{f(re^{i\theta}+c_\nu)}{f(re^{i\theta})}\right|^{l_\nu}
    \end{equation*}
whenever $\theta\in E_1$. Therefore for each $\lambda$ we obtain
    \begin{equation*}
     \int_{E_1} \log^{+}|a_\lambda(re^{i\theta})F_\lambda(re^{i\theta},f)|\,\frac{d\theta}{2\pi} \leq m(r,a_\lambda) + O\left(\sum_{j=1}^{\nu} m\left(r,\frac{f(z+c_j)}{f(z)}\right)\right),
    \end{equation*}
and so, by Theorem \ref{logdiff},
    \begin{equation}\label{E1}
    \int_{E_1} \log^{+}|P(re^{i\theta},f)|\,\frac{d\theta}{2\pi} = o\left(\frac{T(r+|c|,f)^{1+\varepsilon}}{r^\delta}\right) +o(T(r,f))
    \end{equation}
outside of an exceptional set with finite logarithmic measure.

Similarly on $E_2$, by denoting $Q(z,f)=\sum_\gamma b_\gamma(z)
G_\gamma(z,f)$, we obtain
    \begin{equation*}
    \begin{split}
    |P(z,f)| & = \left|\frac{1}{f(z)^n}\sum_\gamma b_\gamma(z) f(z)^{l_0}f(z+c_1)^{l_1}\cdots f(z+c_\mu)^{l_\mu}\right|\\
    &\leq  \sum_\gamma |b_\gamma(z)| \left|\frac{f(re^{i\theta}+c_1)}{f(re^{i\theta})}\right|^{l_1}\cdots\left|\frac{f(re^{i\theta}+c_\mu)}{f(re^{i\theta})}\right|^{l_\mu}
    \end{split}
    \end{equation*}
since $\sum_{j=0}^\mu l_j\leq n$ by assumption. Therefore by Theorem~\ref{logdiff} again,
    \begin{equation}\label{E2}
    \int_{E_2} \log^{+}|P(re^{i\theta},f)|\,\frac{d\theta}{2\pi} = o\left(\frac{T(r+|c|,f)^{1+\varepsilon}}{r^\delta}\right) +o(T(r,f)).
    \end{equation}
The assertion follows by combining \eqref{mP}, \eqref{E1} and
\eqref{E2}.
\end{proof}

Similarly as Theorem \ref{clunieanalogue} can be used to obtain
information about the pole distribution of meromorphic solutions
of difference equations, the next result is concerned with
distribution of \textit{slowly moving targets} $a$ such that
$T(r,a)=o(T(r,f))$ outside of a possible exceptional set of finite
logarithmic measure. In particular, constant functions are always
slowly moving. The following theorem is an analogue of a result
due to A.~Z.~Mohon'ko and V.~D.~Mohon'ko \cite{mohonko:74} on
differential equations.

\begin{theorem}\label{mohonkoanalogue}
Let $f(z)$ be a non-constant meromorphic solution of
    \begin{equation}\label{Peq}
    P(z,f)=0
    \end{equation}
where $P(z,f)$ is difference polynomial in $f(z)$, and let $\delta<1$ and $\varepsilon>0$. If $P(z,a)\not\equiv 0$ for a slowly moving target $a$, then
    \begin{equation*}
    m\left(r,\frac{1}{f-a}\right) = o\left(\frac{T(r+|c|,f)^{1+\varepsilon}}{r^\delta}\right) +o(T(r,f))
    \end{equation*}
for all $r$ outside of a possible exceptional set with finite logarithmic measure.
\end{theorem}

\begin{proof} By substituting $f=g+a$ into \eqref{Peq} we obtain
    \begin{equation}\label{Qeq}
    Q(z,g)+D(z)=0,
    \end{equation}
where $Q(z,g)=\sum_\gamma b_\gamma(z) G_\gamma(z,f)$ is a
difference polynomial in $g$ such that all of its terms are at
least of degree one, and $T(r,D)=o(T(r,g))$ outside a set of
finite logarithmic measure. Also $D\not\equiv0$, since $a$ does not
satisfy \eqref{Peq}. Next we compute $m(r,1/g)$. To this end, note
that the integral to be evaluated vanishes on the part of $|z|=r$
where $|g|>1$. It is therefore sufficient to consider only the
case $|g|\leq 1$. But then,
    \begin{equation*}
    \begin{split}
    \left|\frac{Q(z,g)}{g}\right| & =\frac{1}{|g|} \left| \sum_\gamma b_\gamma(z)g(z)^{l_0}g(z+c_1)^{l_1}\cdots g(z+c_\nu)^{l_\nu}\right|\\
    &\leq  \sum_\gamma |b_\gamma(z)|\left|\frac{g(z+c_1)}{g(z)}\right|^{l_1}\cdots\left|\frac{g(z+c_\nu)}{g(z)}^{l_\nu}\right|\\
    \end{split}
    \end{equation*}
since  $\sum_{j=0}^\nu l_j\geq 1$ for all $\gamma$. Therefore, by equation \eqref{Qeq} and Theorem~\ref{logdiff},
    \begin{equation*}
    \begin{split}
    m\left(r,\frac{1}{g}\right) & \leq m\left(r,\frac{D}{g}\right) + m\left(r,\frac{1}{D}\right)\\
    & = m\left(r,\frac{Q(z,g)}{g}\right) + m\left(r,\frac{1}{D}\right)\\
    & = o\left(\frac{T(r+|c|,g)^{1+\varepsilon}}{r^\delta}\right) +o(T(r,g))
    \end{split}
    \end{equation*}
outside of a set of $r$-values with at most finite logarithmic
measure. Since $g=f-a$ the assertion follows.
\end{proof}

Theorems \ref{clunieanalogue} and \ref{mohonkoanalogue}, like
Theorem \ref{logdiff}, are particularly useful when applied to
functions having finite order. The following two corollaries on
the Nevanlinna deficiency illustrate this fact.

\begin{corollary}\label{cluniecorollary}
Let $f(z)$ be a non-constant finite-order meromorphic solution of
    \begin{equation*}
    f(z)^n P(z,f)=Q(z,f),
    \end{equation*}
where $P(z,f)$ and $Q(z,f)$ are difference polynomials in $f(z)$, and let $\delta<1$.
If the degree of $Q(z,f)$ as a polynomial in $f(z)$ and its shifts
is at most $n$, then
    \begin{equation}\label{mPcor}
    m\big(r,P(z,f)\big) = o\left(\frac{T(r+|c|,f)}{r^\delta}\right) +o(T(r,f))
    \end{equation}
for all $r$ outside of a possible exceptional set with finite logarithmic measure. Moreover,
the Nevanlinna deficiency satisfies
    \begin{equation}\label{dPcor}
    \delta(\infty,P):=\liminf_{r\rightarrow\infty}\frac{m(r,P)}{T(r,P)}= 0.
    \end{equation}
\end{corollary}

\begin{proof} Equation \eqref{mPcor} follows by combining the proof
of Theorem~\ref{clunieanalogue} with Corollary~\ref{logdiffcor},
and so we are left with equation~\eqref{dPcor}. By a well known
result due to Valiron \cite{valiron:31} and A.~Z.~Mohon'ko
\cite{mohonko:71}, we have
    \begin{equation}\label{vmeqn}
    T(r,P)=\deg(P)T(r,f)+o(T(r,f))
    \end{equation}
outside of a possible exceptional set of finite logarithmic measure. In
addition \cite[Lemma 1.1.2]{laine:93} yields
that if $T(r,g)=o(T(r,f))$ outside of an exceptional set of finite
logarithmic measure, then $T(r,g)=o(T(r^{1+\varepsilon},f))$ for any $\varepsilon>0$ and for all $r$ sufficiently
large. Thus, by applying \eqref{mPcor} together
with \eqref{vmeqn} and \cite[Lemma 1.1.2]{laine:93}, we have
    \begin{equation*}
    m(r,P) = o\left(\frac{T(r^{1+\varepsilon},P)}{r^\delta}\right) +o(T(r^{1+\varepsilon},P))
    \end{equation*}
for all sufficiently large $r$. Therefore, since $P$ is of finite
order,
    \begin{equation}\label{m}
    m(r,P)\leq r^{\rho(1+2\varepsilon)-\delta},
    \end{equation}
where $\rho$ is the order of $P$ and $\delta<1$.
Also, there is a sequence $r_n\rightarrow\infty$ as
$n\rightarrow\infty$, such that
    \begin{equation}\label{T}
    T(r_n,P) \geq r_n^{\rho-\varepsilon}
    \end{equation}
for all $r_n$ large enough. The assertion follows by combining
\eqref{m} and \eqref{T} where $\varepsilon$ and $\delta$ are chosen such that $\varepsilon(2\rho+1)<\delta<1$, and by letting $n\rightarrow\infty$.
\end{proof}

\begin{corollary}
Let $f(z)$ be a non-constant finite-order meromorphic solution of
    \begin{equation*}
    P(z,f)=0
    \end{equation*}
where $P(z,f)$ is difference polynomial in $f(z)$, and let $\delta<1$. If $P(z,a)\not\equiv 0$ for a slowly moving target $a$, then
    \begin{equation*}
    m\left(r,\frac{1}{f-a}\right) = o\left(\frac{T(r+|c|,f)}{r^\delta}\right) +o(T(r,f))
    \end{equation*}
for all $r$ outside of a possible exceptional set with finite logarithmic measure. Moreover, the Nevanlinna deficiency satisfies
    \begin{equation*}
    \delta(a,f):=\liminf_{r\rightarrow\infty}\frac{m\left(r,\frac{1}{f-a}\right)}{T(r,f)}=0.
    \end{equation*}
\end{corollary}

We omit the proof since it would be almost identical to that of
Corollary~\ref{cluniecorollary}.

\section{Conclusion}

In this paper we have presented a difference analogue of the Lemma
on the Logarithmic Derivative. This result has potentially a large
number of applications in the study of difference equations. Many
ideas and methods from the theory of differential equations may now
be utilized together with Theorem~\ref{logdiff} to obtain
information about meromorphic solutions of difference equations.
Section~\ref{diffsec} provides a number of examples in this
direction. The analogue of the Clunie Lemma,
Theorem~\ref{clunieanalogue}, may be used to ensure that finite
order meromorphic solutions of certain non-linear difference
equations have a large number of poles. Similarly,
Theorem~\ref{mohonkoanalogue} provides an easy way of telling when a
finite order meromorphic solution of a difference equation does not
have any deficient values.

\end{document}